\documentclass[12pt, english]{article}
\usepackage[T2A]{fontenc}
\usepackage[cp1251]{inputenc}
\usepackage{babel}
\usepackage[centertags]{amsmath}
\usepackage{amsfonts}
\usepackage{amssymb}
\usepackage{amsthm}
\usepackage{newlfont}
\usepackage[active]{srcltx}

\newtheorem{theorem}{\bf Theorem}
\newtheorem{remark}{\bf Remark}

\begin{document}

\title{
Esseen type bounds of the remainder in a combinatorial CLT }

\author{Andrei N. Frolov
\\ Dept. of Mathematics and Mechanics
\\ St.~Petersburg State University
\\ St. Petersburg, Russia
\\ E-mail address: Andrei.Frolov@pobox.spbu.ru}

\maketitle

{\abstract{ We derive 
Esseen type bounds of the remainder in a combinatorial central limit theorem
for independent random variables 
without third moments.
 }


}

\section{Introduction and results}


Let $\|c_{ij}\|$ be a $n \times n$ matrix of real numbers such that
\begin{eqnarray}\label{10}
\sum\limits_{i=1}^n c_{ij} = \sum\limits_{j=1}^n c_{ij}=0.
\end{eqnarray}

Let $\|Y_{ij}\|$ be a $n \times n$ matrix of independent random variables
with $E Y_{ij}=0$ and $E Y_{ij}^2=\sigma_{ij}^2$.

Assume that $\pi=(\pi(1),\pi(2),\ldots\pi(n))$ is a random permutation of
$1,2,\ldots, n$ with the uniform distribution on the set of all
such permutations.
Assume that $\pi$ and $\|Y_{ij}\|$ are independent.

Denote
$$ S_n = \sum\limits_{i=1}^n c_{i\pi(i)} + \sum\limits_{i=1}^n Y_{i\pi(i)}.
$$
Suppose that $n\geqslant 2$. Then
$$ E S_n=0, \quad
B_n= D S_n = \frac{1}{n-1}\sum\limits_{i,j=1}^n c_{ij}^2 +
\frac{1}{n} \sum\limits_{i,j=1}^n \sigma_{ij}^2.
$$
To avoid triviality, we assume in the sequel that $B_n>0$.

Put $$  \Delta_n = \sup\limits_{x\in \mathbb{R}} \left| P(S_n<
x\sqrt{B_n})-\Phi(x)\right|, $$
where $\Phi(x)$ is the standard normal distribution function.

First results on asympotical normality of $S_n$ were obtained
in the case $P(Y_{ij}=0)=1$ for $1\leqslant i, j \leqslant n$.
Motivated by statistical applications, Wald and Wolfowits (1944)
stated conditions, sufficient for $\Delta_n \rightarrow 0$
as $n \rightarrow\infty$ when
$c_{ij}=a_i b_j$. Noether (1949) proved that these conditions maybe
replaced by weaker ones. Hoeffding (1951) considered general case
of $c_{ij}$ and obtained a combinatorial central limit theorem (CLT).
Further results on combinatorial CLT were obtained by Motoo (1957)
and Kolchin and Chistyakov (1973).
Von Bahr (1976) and Ho and Chen (1978) derived bounds for the remainder
in combinatorial CLT in the case on non-degenerated $Y_{ij}$.
Botlthausen (1984) obtained Esseen type inequality for the remainder in
the case of degenerated $Y_{ij}$. The constant was not be specified
in the latter paper. Further results of this type were proved
by Goldstein (2005) (see also Chen, Goldstein and Shao (2011)).
They contain explicit constants in the inequalities.
For non-degenerated $Y_{ij}$,  Esseen type inequalities were
stated by Neammanee and Suntornchost (2005), Neammanee and Rattanawong (2009)
 and Chen and Fang (2012) (see also comments on p.2 of Chen and Fang (2012)).

At the moment, best results on bounds in combinatorial CLT are obtained
by Stein's method. The detailed discussion of this approach may be
found in Chen and Fang (2012) and references therein.

We start with the following known result.

\medskip \noindent
{\bf Theorem A.}
{\it If $E |Y_{ij}|^3 <\infty$ for $1\leqslant i, j \leqslant n$ then
there exists an absolute positive constant $A_0$ such that
\begin{equation}\label{be}
  \Delta_n \leqslant
\frac{A_0}{ B_n^{3/2} n} \sum\limits_{i,j=1}^n |c_{ij}+Y_{ij}|^3.
\end{equation}
}
\medskip

Chen and Fang (2012) have proved that inequality (\ref{be}) holds with $A_0=447$.

In this paper, we derive generalizations of (\ref{be}) to the case of
random variables $Y_{ij}$ without third moments. We apply the classical technique
of truncation in a similar way as in Petrov (1995) where generalizations of
the Esseen's inequality may be found.

Our first result is as follows.

\begin{theorem}\label{th1}
{\it There exists an absolute positive constant $A$ such that
\begin{equation}\label{be1}
  \Delta_n \leqslant A \left(
\frac{1}{ B_n^{3/2} n} \sum\limits_{i,j=1}^n |c_{ij}|^3 +
\frac{1}{B_n n} \sum\limits_{i,j=1}^n \alpha_{ij}
+ \frac{1}{ B_n^{3/2}  n} \sum\limits_{i,j=1}^n \beta_{ij}
\right),
\end{equation}
where $\alpha_{ij}=E Y_{ij}^2 I\{|Y_{ij}| \geqslant \sqrt{B_n}\}$ and
$\beta_{ij}=E\left|Y_{ij} \right|^3 I\{|Y_{ij}| < \sqrt{B_n} \}$
for $1\leqslant i, j \leqslant n$,
and $I\{\cdot\}$ denotes the indicator of the event in brackets.
}
\end{theorem}

\begin{remark}\label{r1}
One can put $A=\max\left\{1810, 198 A_0+5\right\}$ in Theorem \ref{th1}.
\end{remark}

Theorem \ref{th1} implies the next result.

\begin{theorem}\label{th2}
{\it Let $g(x)$ be a positive, even function such that
$g(x)$ and $x/g(x)$ are non-decreasing for $x>0$. Suppose that
$g_{ij}=E Y_{ij}^2 g(Y_{ij})<\infty$ for $1\leqslant i, j \leqslant n$.
Then
\begin{equation}\label{be2}
  \Delta_n 
\leqslant 2 A \left( \frac{1}{ B_n^{3/2} n} \sum\limits_{i,j=1}^n |c_{ij}|^3 +
\frac{1}{B_n g(\sqrt{B_n}) n} \sum\limits_{i,j=1}^n g_{ij} 
\right).
\end{equation}
}
\end{theorem}

It is clear that
$\alpha_{ij} \leqslant g_{ij}/g(\sqrt{B_n})$ and
$\beta_{ij} \leqslant \sqrt{B_n} g_{ij}/g(\sqrt{B_n})$
for $1\leqslant i, j \leqslant n$.
Hence (\ref{be1}) implies (\ref{be2}). At the same time, putting
$g(x)=\min\{|x|,\sqrt{B_n}\}$ in Theorem \ref{th2}, we conclude
that (\ref{be2}) yields (\ref{be1}) with $2 A$ instead of $A$.

One can obtain variants of Theorem \ref{th2} as follows. Writing
$$ \alpha_{ij}=E Y_{ij}^2 I\{Y_{ij} \leqslant -\sqrt{B_n}\}+
E Y_{ij}^2 I\{Y_{ij} \geqslant \sqrt{B_n}\}
$$
allows estimate $\alpha_{ij}$ by a sum of different moments of
random variables $Y_{ij}^-$ and $Y_{ij}^+$, where
$Y_{ij}^\pm=\max\{\pm Y_{ij}, 0\}$. The same maybe done
also for $\beta_{ij}$. So, if
$g^-(x)$ and $g^+(x)$, $x>0$, are positive, non-decreasing functions such that
$x/g^-(x)$ and $x/g^+(x)$ are non-decreasing, then
$$ \alpha_{ij} \leqslant \frac{g_{ij}^-}{g^-(\sqrt{B_n})}+
\frac{g_{ij}^+}{g^+(\sqrt{B_n})} \quad \mbox{and}
\quad
\beta_{ij} \leqslant \frac{\sqrt{B_n} g_{ij}^-}{g^-(\sqrt{B_n})}+
\frac{\sqrt{B_n} g_{ij}^+}{g^+(\sqrt{B_n})},
$$
where
$g_{ij}^- = E (Y_{ij}^-)^2 g^-(Y_{ij}^-)$ and
$g_{ij}^+ = E (Y_{ij}^+)^2 g^+(Y_{ij}^+)$.
Further generalizations maybe derived in the same way for
$g^-(x)$ and $g^+(x)$ depending on $i$ and $j$.

One of the most important case of $g(x)$ is $g(x)=|x|^\delta$,
$\delta \in (0,1]$, in which Theorem \ref{th2} turns to the
following result.

\begin{theorem}\label{th3}
{\it Assume that
$E |Y_{ij}|^{2+\delta}<\infty$ for $1\leqslant i, j \leqslant n$,
where $\delta \in (0,1]$. Then
\begin{equation}\label{be3}
  \Delta_n \leqslant
  2 A \left( \frac{1}{ B_n^{3/2} n} \sum\limits_{i,j=1}^n |c_{ij}|^3 +
\frac{1}{B_n^{1+\delta/2}  n} \sum\limits_{i,j=1}^n E |Y_{ij}|^{2+\delta}
\right).
\end{equation}
}
\end{theorem}

Generalizations of Theorem \ref{th3} maybe obtained in the same way as
it was mentioned for Theorem \ref{th2} above.

Inequality (\ref{be}) is better than (\ref{be3}) for $\delta=1$. But if
$P(Y_{ij}=0)=1$ for $1\leqslant i, j \leqslant n$, then right-hand sides
of these bounds coincide up to constants.

The second term in righthand side of (\ref{be3}) is the Lyapunov type ratio
for $Y_{ij}$. Theorem \ref{th3} yields combinatorial CLT under Lyapunov
type condition on $Y_{ij}$. Theorem \ref{th1} also implies combinatorial CLT
under Lindeberg type condition on $Y_{ij}$. The latter follows from
(\ref{be1}) and
$$ \beta_{ij} \leqslant \sqrt{B_n}
\left( E Y_{ij}^2 I\{|Y_{ij}| \geqslant \varepsilon \sqrt{B_n}\}
+ \varepsilon \sigma_{ij}^2 \right)
$$
for $1\leqslant i, j \leqslant n$ and all $\varepsilon \in (0,1)$.

\section{Proofs}

{\bf Proof of Theorem \ref{th1}.}
Assume first that $B_n=1$.

Put
\begin{eqnarray*}
&&
\bar{Y}_{ij} = Y_{ij} I\{|Y_{ij}| < 1 \},\quad
\bar{a}_{ij} = E \bar{Y}_{ij},\quad \bar{\sigma}^2_{ij} = D \bar{Y}_{ij},
\\ &&
\bar{a}_{i.} =\frac{1}{n} \sum\limits_{j=1}^n \bar{a}_{ij},
\quad
\bar{a}_{.j} =\frac{1}{n} \sum\limits_{i=1}^n \bar{a}_{ij},
\quad
\bar{a}_{..} =\frac{1}{n^2} \sum\limits_{i,j=1}^n \bar{a}_{ij},
\end{eqnarray*}
for $1\leqslant i, j \leqslant n$.
Denote
$$
\bar{S}_n = \sum\limits_{i=1}^n c_{i\pi(i)}
+ \sum\limits_{i=1}^n \bar{Y}_{i\pi(i)}.
$$
Hence
$$ \bar{e}_n=E \bar{S}_n= n \bar{a}_{..}, \quad
\bar{B}_n= D \bar{S}_n = \frac{1}{n-1}\sum\limits_{i,j=1}^n
(c_{ij}+\bar{a}_{ij}-\bar{a}_{i.}-\bar{a}_{.j}+\bar{a}_{..})^2
+\frac{1}{n} \sum\limits_{i,j=1}^n \bar{\sigma}_{ij}^2.
$$
Put
\begin{eqnarray*}
F_n(x)=P\left(S_n < x \right), \quad
\bar{F}_n(x)=P\left(\bar{S}_n < x \right).
\end{eqnarray*}
Denote
\begin{eqnarray*}
&&
\Delta_{n1} = \sup\limits_{x\in \mathbb{R}} |F_n(x)-\bar{F}_n(x)|,\quad
\Delta_{n2} = \sup\limits_{x\in \mathbb{R}}
|\bar{F}_n(x)-\Phi((x-\bar{e}_n)/\sqrt{\bar{B}_n})|, 
\\ &&
\Delta_{n3} = \sup\limits_{x\in \mathbb{R}}
|\Phi((x-\bar{e}_n)/\sqrt{\bar{B}_n}) - \Phi(x/\sqrt{\bar{B}_n})|,
\quad
\Delta_{n4} = \sup\limits_{x\in \mathbb{R}} |\Phi(x/\sqrt{\bar{B}_n})-\Phi(x)|.
\end{eqnarray*}
It is clear that
$$ \Delta_n= \sup\limits_{x\in \mathbb{R}} |F_n(x)-\Phi(x)|
\leqslant \sum\limits_{k=1}^4 \Delta_{nk}.
$$

Denote
$$
C_n= \frac{1}{n} \sum\limits_{i,j=1}^n |c_{ij}|^3, \quad
L_n= \frac{1}{n} \sum\limits_{i,j=1}^n
E Y_{ij}^2 I\{|Y_{ij}| \geqslant 1\}.
$$

We have
\begin{eqnarray*}
&&
\hspace*{-\parindent}
1-\bar{B}_n=
B_n-\bar{B}_n
\\ &&
\hspace*{-\parindent}
= \frac{1}{n } \sum\limits_{i,j=1}^n ( \sigma_{ij}^2 - \bar{\sigma}_{ij}^2) +
\frac{1}{n-1 } \sum\limits_{i,j=1}^n
(c_{ij}^2-
(c_{ij}+\bar{a}_{ij}-\bar{a}_{i.}-\bar{a}_{.j}+\bar{a}_{..})^2)
\\ &&
\hspace*{-\parindent}
=\frac{1}{n } \sum\limits_{i,j=1}^n ( \sigma_{ij}^2 - \bar{\sigma}_{ij}^2) -
\frac{1}{n-1 } \sum\limits_{i,j=1}^n
(\bar{a}_{ij}-\bar{a}_{i.}-\bar{a}_{.j}+\bar{a}_{..})^2
- \frac{2}{n-1 } \sum\limits_{i,j=1}^n
c_{ij} (\bar{a}_{ij}-\bar{a}_{i.}-\bar{a}_{.j}+\bar{a}_{..})
\\ &&
\hspace*{-\parindent}
=\frac{1}{n } \sum\limits_{i,j=1}^n ( \sigma_{ij}^2 - \bar{\sigma}_{ij}^2) -
\frac{1}{n-1 } \sum\limits_{i,j=1}^n
(\bar{a}_{ij}-\bar{a}_{i.}-\bar{a}_{.j}+\bar{a}_{..})^2
- \frac{2}{n-1 } \sum\limits_{i,j=1}^n
c_{ij} \bar{a}_{ij}.
\end{eqnarray*}
In the last equality, we have used (\ref{10}).

In the sequel, we will repeatedly use that
$|\bar{a}_{ij}|< 1$ for $1\leqslant i, j \leqslant n$.

We have
\begin{eqnarray}\label{20}
0 \leqslant \sigma_{ij}^2 - \bar{\sigma}_{ij}^2 =
E Y_{ij}^2 I\{|Y_{ij}| \geqslant 1\} + (\bar{a}_{ij})^2 \leqslant
E Y_{ij}^2 I\{|Y_{ij}| \geqslant 1\} + |\bar{a}_{ij}|
\end{eqnarray}
for $1\leqslant i, j \leqslant n$.

Taking into account that $(x+y)^2 \leqslant 2 (x^2+y^2)$
for all real $x$ and $y$, we get
\begin{eqnarray}\label{30}
 (\bar{a}_{ij}-\bar{a}_{i.}-\bar{a}_{.j}+\bar{a}_{..})^2
\leqslant 4 (\bar{a}_{ij}^2+\bar{a}_{i.}^2+\bar{a}_{.j}^2+\bar{a}_{..}^2)
\leqslant 4 (|\bar{a}_{ij}|+|\bar{a}_{i.}|+|\bar{a}_{.j}|+|\bar{a}_{..}|)
\end{eqnarray}
for $1\leqslant i, j \leqslant n$.
Moreover,
\begin{eqnarray}\label{40}
\sum\limits_{i,j=1}^n |\bar{a}_{i.}| \leqslant
\sum\limits_{i,j=1}^n \left(\frac{1}{n} \sum\limits_{j=1}^n |\bar{a}_{ij}|\right)
= \sum\limits_{i,j=1}^n |\bar{a}_{ij}|
\end{eqnarray}
and, in the same way,
\begin{eqnarray}\label{50}
\sum\limits_{i,j=1}^n |\bar{a}_{.j}| \leqslant
 \sum\limits_{i,j=1}^n |\bar{a}_{ij}|.
\end{eqnarray}
Further,
\begin{eqnarray}\label{60}
\sum\limits_{i,j=1}^n |\bar{a}_{..}|\leqslant
\sum\limits_{i,j=1}^n
\left( \frac{1}{n^2} \sum\limits_{i,j=1}^n |\bar{a}_{ij}|\right)
= \sum\limits_{i,j=1}^n |\bar{a}_{ij}|.
\end{eqnarray}

Noting that $ x y \leqslant x^{3}/3 + 2 y^{3/2}/3$ for all non-negative $x$ and $y$,
we get
\begin{eqnarray}\label{70}
\sum\limits_{i,j=1}^n |c_{ij} \bar{a}_{ij}|
\leqslant \frac{1}{3} \sum\limits_{i,j=1}^n |c_{ij}|^3+
\frac{2}{3} \sum\limits_{i,j=1}^n |\bar{a}_{ij}|^{3/2}
\leqslant \frac{1}{3} \sum\limits_{i,j=1}^n |c_{ij}|^3+
\frac{2}{3} \sum\limits_{i,j=1}^n |\bar{a}_{ij}|.
\end{eqnarray}

Making use of $E Y_{ij}=0$,
\begin{eqnarray} \label{80}
 |\bar{a}_{ij}| = |E Y_{ij} I\{|Y_{ij}| < 1 \}|
 = |E Y_{ij} I\{|Y_{ij}| \geqslant 1 \}| \leqslant
 E |Y_{ij}| I\{|Y_{ij}| \geqslant 1 \}
 \leqslant  E Y_{ij}^2 I\{|Y_{ij}| \geqslant 1 \}
\end{eqnarray}
for $1\leqslant i, j \leqslant n$.

It follows by (\ref{20})--(\ref{80}) that
\begin{eqnarray*}
&& | B_n-\bar{B}_n| \leqslant
\left(2+ \frac{16 n}{n-1} + \frac{2 n}{3(n-1)} \right) L_n
+ \frac{n}{3(n-1)} C_n.
\end{eqnarray*}
This yields that for $n\geqslant 2$,
\begin{eqnarray}\label{90}
 | B_n-\bar{B}_n| \leqslant 36 L_n + C_n.
\end{eqnarray}

Assume that $\bar{B}_n \leqslant \varrho B_n$ where $\varrho \in (0,1).$
Then
\begin{equation}\label{95}
  \frac{36 L_n + C_n}{ 1-\varrho}
 \geqslant 1 \geqslant \Delta_n
\end{equation}
and the conclusion of Theorem \ref{th1} holds.

In the sequel, we assume that $\bar{B}_n > \varrho B_n$.

Now we will estimate $\Delta_{nk}$, $k=\overline{1,4}.$

Since
$$ \left\{S_n < x \right\} \subset
\left\{\bar{S}_n < x \right\} \cup
\bigcup\limits_{i=1}^n \left\{|Y_{i\pi(i)}| \geqslant
1 \right\},
$$
we have
\begin{eqnarray*}
&& F_n(x) \leqslant \bar{F}_n(x)+ \sum\limits_{i=1}^n
P\left(|Y_{i\pi(i)}| \geqslant 1 \right)
= \bar{F}_n(x)+ \frac{1}{n} \sum\limits_{i,j=1}^n
P\left(|Y_{ij}| \geqslant 1 \right)
\\ &&
\leqslant \bar{F}_n(x)+ \frac{1}{n}
\sum\limits_{i,j=1}^n  E Y_{ij}^2 I\{|Y_{ij}| \geqslant 1 \}
= \bar{F}_n(x)+L_n.
\end{eqnarray*}
From the other hand
$$ \left\{\bar{S}_n < x \right\} \subset
\left\{S_n < x \right\} \cup
\bigcup\limits_{i=1}^n \left\{|Y_{i\pi(i)}| \geqslant
1 \right\},
$$
which yields that
$$ \bar{F}_n(x) \leqslant F_n(x)+ L_n.
$$
It follows that
\begin{eqnarray}\label{100}
 \Delta_{n1} \leqslant  L_n.
\end{eqnarray}

Applying Theorem A with
$c_{ij}+\bar{a}_{ij}-\bar{a}_{i.}-\bar{a}_{.j}+\bar{a}_{..}$ and
$\bar{Y}_{ij}-\bar{a}_{ij}$ instead of
$c_{ij}$ and $Y_{ij}$ correspondingly, we get inequality
\begin{eqnarray*}
&& \Delta_{n2} \leqslant A_0 \frac{1}{n \bar{B}_n^{3/2}}
\sum\limits_{i,j=1}^n E\left|c_{ij}+ \bar{Y}_{ij}
-\bar{a}_{i.}-\bar{a}_{.j}+\bar{a}_{..} \right|^3.
\end{eqnarray*}
Using that $|x+y|^3 \leqslant 4(|x|^3+|y|^3)$ for all real $x$ and $y$,
we obtain
\begin{eqnarray}
&& \nonumber
\hspace*{-\parindent}
\Delta_{n2}
\leqslant  A_0 \frac{64}{n \varrho^{3/2}}
\sum\limits_{i,j=1}^n
(|c_{ij}|^3+ E\left|\bar{Y}_{ij} \right|^3+
|\bar{a}_{i.}|^3 +|\bar{a}_{.j}|^3+|\bar{a}_{..}|^3
)
\\ &&
\hspace*{-\parindent}
\label{110}
\leqslant  A_0 \frac{64}{n \varrho^{3/2}}
\sum\limits_{i,j=1}^n
(|c_{ij}|^3+ E\left|\bar{Y}_{ij} \right|^3+
|\bar{a}_{i.}| +|\bar{a}_{.j}|+|\bar{a}_{..}|
)
\leqslant  A_0 \frac{64}{\varrho^{3/2}} ( C_n + 3 L_n+D_n),
\end{eqnarray}
where
$$ D_n = \frac{1}{n} \sum\limits_{i,j=1}^n  E\left|\bar{Y}_{ij} \right|^3 =
\frac{1}{n} \sum\limits_{i,j=1}^n
E\left|Y_{ij} \right|^3 I\{|Y_{ij}| < 1 \}.
$$

The inequality
$$ \sup\limits_{x\in \mathbb{R}} |\Phi(x+y)-\Phi(x)|
\leqslant \frac{|y|}{\sqrt{2 \pi }}
$$
holds for all real $y$.
By the latter inequality and (\ref{80}), we have
\begin{eqnarray}\label{120}
\Delta_{n3} \leqslant  \frac{|\bar{e}_n|}{\sqrt{2 \pi } \bar{B}_n^{1/2}}
\leqslant \frac{|\bar{e}_n|}{\sqrt{2 \pi } \varrho^{1/2}} \leqslant
\frac{1}{\sqrt{2 \pi } \varrho^{1/2} n}
\sum\limits_{i,j=1}^n |\bar{a}_{ij}| \leqslant
\frac{ L_n}{\sqrt{2 \pi } \varrho^{1/2} }.
\end{eqnarray}

Note that
$$ \sup\limits_{x\in \mathbb{R}} |\Phi(xy)-\Phi(x)|
\leqslant \frac{y-1}{\sqrt{2 \pi } e}
$$
for all $y \geqslant 1$.
If $\bar{B}_n < B_n=1$, then
$$ \Delta_{n4} \leqslant
\frac{1}{\sqrt{2 \pi e}} \left(\frac{1}{\sqrt{\bar{B}_n}}-1\right)
$$
and we conclude by (\ref{90}) that
$$ \frac{1}{\sqrt{\bar{B}_n}}-1
=\frac{B_n-\bar{B}_n}{\sqrt{\bar{B}_n}(\sqrt{B_n}
+\sqrt{\bar{B}_n})}
\leqslant \frac{36 L_n +C_n}{\varrho^{1/2}(1+\varrho^{1/2})}.
$$

If $y\in(0,1)$, then
$$ \sup\limits_{x\in \mathbb{R}} |\Phi(xy)-\Phi(x)|
\leqslant \frac{1-y}{y \sqrt{2 \pi } e}.
$$
So, for $\bar{B}_n \geqslant B_n=1$, inequality
$$ \Delta_{n4} \leqslant
\frac{1}{\sqrt{2 \pi e}} \left(\sqrt{\bar{B}_n}-1\right)
$$
holds and we have
$$ \sqrt{\bar{B}_n}-1 = \frac{B_n-\bar{B}_n}{\sqrt{B_n}
+\sqrt{\bar{B}_n}}.
$$

It follows that
\begin{eqnarray}\label{130}
 \Delta_{n4} \leqslant
\frac{36 L_n +C_n}{\sqrt{2 \pi e} \varrho^{1/2}(1+\varrho^{1/2})}.
\end{eqnarray}

Relations (\ref{100})--(\ref{130}) imply
$$ \Delta_n \leqslant A_1 (L_n+C_n+D_n),
$$
when $B_n=1$. If $B_n\neq 1$, we replace $c_{ij}$ and $Y_{ij}$ by
$c_{ij}/\sqrt{B_n}$ and $Y_{ij}/\sqrt{B_n}$ in the previous part of the proof.
Then we obtain from the last inequality that
$$ \Delta_n \leqslant A_1 \left(
\frac{1}{ n} \sum\limits_{i,j=1}^n \frac{|c_{ij}|^3}{ B_n^{3/2}} +
\frac{1}{ n} \sum\limits_{i,j=1}^n
E \frac{Y_{ij}^2}{B_n} I\left\{\frac{|Y_{ij}|}{\sqrt{B_n}} \geqslant 1\right\}
+ \frac{1}{  n } \sum\limits_{i,j=1}^n
E\frac{\left|Y_{ij} \right|^3}{B_n^{3/2}} I\left\{\frac{|Y_{ij}|}{\sqrt{B_n}} <1  \right\}
\right),
$$
which coincides with (\ref{be1}).
$\Box$

\medskip
{\bf Proof of Remark \ref{r1}.}
From (\ref{95})--(\ref{130}) we get
$$ A_1 = \max\left\{\frac{36}{1-\varrho},\,
1+ A_0 \frac{192}{\varrho^{3/2}}+\frac{1}{\sqrt{2\pi} \varrho^{1/2}}
+ \frac{36}{\sqrt{2\pi} e \varrho^{1/2} (1+\varrho^{1/2})}
\right\}.
$$
This yields that
$$ A_1 < \max\left\{\frac{36}{1-\varrho},\,
1+ A_0 \frac{192}{\varrho^{3/2}}+\frac{0.4}{ \varrho^{1/2}}
+ \frac{5.3}{ \varrho^{1/2} (1+\varrho^{1/2})}
\right\}.
$$
Taking $\varrho^{1/2}=0.99$, we get $A_1 <\max\left\{1810, 198 A_0+5\right\}.$
$\Box$

\bigskip
\noindent
{\bf References}
{\footnotesize

\parindent 0 mm

Bolthausen E. (1984) An estimate of the remainder in
a combinatorial central limit theorem.
Z. Wahrsch. verw. Geb. 66, 379-386.

Goldstein L. (2005) Berry-Esseen bounds for combinatorial central limit theorems
and pattern occurrences, using zero and size biasing. J. Appl. Probab. 42, 661-683.

Chen L.H.Y., Fang X. (2012) 0n the error bound in a
combinatorial central limit theorem. arXiv:1111.3159.

Chen L.H.Y., Goldstein L., Shao Q.M. (2011) Normal approximation
by Stein's method. Springer.

Ho S.T., Chen L.H.Y. (1978) An $L_p$ bounds for the remainder in
a combinatorial central limit theorem.
Ann. Probab. 6, 231-249.

Hoeffding W. (1951) A combinatorial central limit theorem. Ann. Math. Statist. 22,
558-566.

Kolchin V.F., Chistyakov V.P. (1973) On a combinatorial central limit theorem.
Theor. Probab. Appl. 18, 728-739.

Motoo M. (1957) On Hoeffding's combinatorial central limit theorem.
Ann. Inst. Statist. Math. 8, 145-154.

Neammanee K., Suntornchost J. (2005) A uniform bound on
a combinatorial central limit theorem.
Stoch. Anal. Appl. 3, 559-578.

Neammanee K., Rattanawong P. (2009) A constant on a uniform bound of
a combinatorial central limit theorem. J. Math. Research 1, 91-103.

Noether G.E. (1949) On a theorem by Wald and  Wolfowitz.
 Ann. Math. Statist. 20, 455-458.

Petrov V.V. (1995) Limit theorems of probability theory. Sequences of
independent random variables. Clarendon press, Oxford.

von Bahr B. (1976) Remainder term estimate in a
combinatorial central limit theorem.
Z. Wahrsch. verw. Geb. 35, 131-139.

Wald A., Wolfowitz J. (1944) Statistical tests based on permutations of
observations. Ann. Math. Statist. 15, 358-372.

}

\end{document}